\documentclass[10pt,letterpaper]{article}
\usepackage{blindtext,titlefoot}
\usepackage{authblk}
\usepackage{marvosym}
\usepackage[utf8]{inputenc}
\usepackage[english]{babel}
\usepackage{amsmath}
\usepackage{amsfonts}
\usepackage{amssymb}
\usepackage{makeidx}
\usepackage{graphicx}
\usepackage{lmodern}
\usepackage{kpfonts}
\usepackage{dsfont}
\usepackage{cancel}
\usepackage{amsthm} 
\usepackage[left=2cm,right=2cm,top=2cm,bottom=2cm]{geometry}
\usepackage{subcaption}
\usepackage{multirow}
\usepackage{float}
\usepackage{url}
\usepackage{breakurl}
\usepackage[breaklinks]{hyperref}
\usepackage{multicol}
\usepackage{mathtools, cuted}
\usepackage{lipsum}
\usepackage{booktabs}
\usepackage{comment}
\usepackage{cite}
\usepackage{fancyhdr}
\fancypagestyle{firstpage}{\fancyhf{}
\setcounter{page}{1}
\rfoot{\thepage}
}

\providecommand{\set}[1]{
\left\{#1\right\}
}

\usepackage{enumitem} % Customized lists
\setlist[itemize]{noitemsep} % Make itemize lists more compact

\usepackage{titlesec} % Allows customization of titles
\titleformat{\section}[block]{\large\bfseries\scshape\centering}{\thesection.}{1em}{} % Change the look of the section titles
\titleformat{\subsection}[block]{\large\bfseries\scshape\centering}{\thesubsection.}{1em}{}
\titleformat{\subsubsection}[block]{\large\bfseries\scshape\centering}{\thesubsubsection.}{1em}{} % Change the look of the section titles

\title{\huge\bfseries  Fractional flow equations.\\ A model for pressure deficit in an oil well.
}

\author[,a]{\normalsize B. F. Martínez-Salgado  \footnote{E-mail: masabemx@yahoo.com.mx}}

\author[,a]{\normalsize F. Alcántara-López  \footnote{E-mail: alcantaralopezfernando@gmail.com}}

\author[,b]{\normalsize A. Torres-Hernandez  \footnote{E-mail: anthony.torres@ciencias.unam.mx}}

\author[,a]{\\ \normalsize F. Brambila-Paz \footnote{E-mail: fernandobrambila@gmail.com}}

\author[,c]{\normalsize C. Fuentes \footnote{E-mail: cbfuentesr@gmail.com}}

\author[,a]{\normalsize J. López Estrada  \footnote{E-mail: jelpze@gmail.com}}

\affil[a]{\normalsize Department of Mathematics, Faculty of Science - UNAM, Mexico}

\affil[b]{\normalsize Department of Physics, Faculty of Science - UNAM, Mexico}

\affil[c]{\normalsize Mexican Institute of Water Technology, Jiutepec, Morelos, Mexico}

\date{}

\begin{document}

% Print the title
\maketitle

\thispagestyle{firstpage}

%\maketitle\unmarkedfntext{DOI :10.5121/mathsj.2020.7102}
%----------------------------------------------------------------------------------------
%	ARTICLE CONTENTS
%----------------------------------------------------------------------------------------

\begin{abstract}

This article presents a novel system of flow equations that models the pressure deficit of a reservoir considered as a triple continuous medium formed by the rock matrix, vugular medium and fracture. In non-conventional reservoirs, the velocity of the fluid particles is altered due to physical and chemical phenomena caused by the interaction of the fluid with the medium, this behavior is defined as anomalous. A more exact model can be obtained with the inclusion of the memory formalism concept that can be expressed through the use of fractional derivatives. Using Laplace transform of the Caputo fractional derivative and Bessel functions, a semi-analytical solution is reached in the Laplace space.

\end{abstract}

\section*{Introduction} % This section will not appear in the table of contents due to the star (\section*)
The  modeling of an oil reservoir is of paramount importance, since it allows to take decisions   that can improve the extraction of hydrocarbons. Over the years, various approaches have been used to obtain a more comprehensive model of the fluid behavior within the reservoir through the interpretation of pressure deficit.
Here is presents a new semi-analytical solution of a system of partial fractional differential equations of flows coupled in a system with triple porosity and triple permeability using Caputo type temporary fractional derivatives. This system is reduced by Laplace transforms to one where  each equation can be see as a Bessel type of the second type. This system is solved in the Laplace space using algebraic method, the necessary definitions and concepts are given previously such that fractional calculus, Bessel equations and Laplace transforms.
\newline
%\keywords{Caputo derivative, Laplace Transforms, Bessel equations}
\section{Previous Works}
As mentioned above, many approximations have been developed to model correctly the fluid behavior in an oil reservoir. Warren \& Root in \cite{Warren1963} proposed equations where they considered that  matrix and the fractures systems had an Euclidean structure.
\newline

 From this approach, Chang and Yortsos in \cite{Chang1990} presented a formulation where a fractal fracture in a Euclidean matrix is considered. Camacho-Velazquez et al. in \cite{CamachoVelazquez2008} take up this issue again,  and proposed a model of double porosity in naturally fractured vugular sites. In their model they make use of a fractional order Caputo-type derivative, which has already been proposed in flow models by Metzler-Glökle-Nonenmacher as can be seen in \cite{Metzler1994}. Camacho et al. in their article \cite{CamachoVelazquez2005}  generalized a classical flow equation to  one that considers the medium as the union of two or three porous media (fractured ,  vugular  and  matrix media for the latter case).
 \newline
 
  Classical models are constructed from the principle of conservation of the mass of each of the fluids involved in the same media and the Darcy's law for fluid in porous media as illustrated by Ertekin in \cite{Ertekin1984}. Martínez Salgado et al \cite{martinez2017applications2}  developed a model for a triadic medium with triple porosity and triple permeability using Caputo fractional derivatives for time and Weyl fractional derivatives for space. In \cite{fernando2017fractional,brambila2018fractional
  ,torreshern2020,torres2020reduction,torres2020fractional}  fractional derivatives were used in the search of
  complex roots of non-linear equations by fractional iterative methods. Martínez et al solved in \cite{martineznumerical} and \cite{martinez2017applications1} numerically diffusion equations with Riesz fractional derivative  in space by radial-based functions.  The equation with fractional derivative uses a fractional Darcy's law deduced by Le Mehaute as seen in \cite{LeMehaute1984} in the same way that appears in Raghavan's article, \cite{Raghavan2013}, where the order of fractional derivatives is expressed in terms of the Hausdorff dimension of the medium. Furthermore, taking the fractional Darcy Law, the order can be obtained from data  \cite{Holy2016} or through inverse problems  \cite{hatano2013determination}.\\

\section{Methods}
The classical model assumes that the properties of rock and fluids are stable, the hydrodynamics of the flow of fluids in the porous medium is adequately described by Darcy's law, the geometry of the reservoir is of the Euclidean type. 

The basis of the model is found in the continuity equation and Darcy's law for a flow through a porous medium, as illustrated in \cite{Bear1988} and  \cite{Peaceman1977}, these equations can be expressed as

\begin{align}
&\dfrac{\partial (\rho \theta)}{\partial t}+\nabla\cdot p(\rho q)=\rho\Upsilon,\nonumber \\ 
& q=-\dfrac{1}{\mu}k(p)(\nabla p-\rho g\nabla D),
\end{align}
where $\theta$  is the volumetric content of the fluid; $ q = (q_ {1}, q_ {2}, q_ {3}) $ is the flow of Darcy; with its spatial components $ (x, y, z ) $, $ t $ is the time; $ \rho$ is the density of the fluid; $ \mu $ is the dynamic viscosity of the fluid; $ g $  is the gravitational acceleration, $ \Upsilon $ is a source term and represents a volume contributed by fluid per unit volume of porous medium in the unit of time; $ p $ is the pressure; $ D $ is the depth as a function of spatial coordinates, generally assimilated to the vertical coordinate $ z $; $ k $ is the permeability tensor of the porous medium, $ \theta(p) $ and $ k (p) $ are characteristics of the fluid dynamics of the medium.

The general equation of fluid transfer is obtained by combining the equations as in  \cite{Fuentes2014}:
\begin{equation*}
\dfrac{\partial(\rho\theta)}{\partial t}=\nabla\cdot p[\frac{\rho}{\mu}k(p)(\nabla p-\rho g\nabla D)]+\rho\Upsilon. 
\end{equation*}

This differential equation contains two dependent variables, namely the moisture content $ \theta $ and the fluid pressure $ p $, wich are related. For this reason, the saturation $ S (p) $ is defined as
\begin{equation}
\theta(p)=\phi(p)S(p),
\end{equation}
where $ \phi $ is the total porosity of the medium. The specific capacity is defined by
\begin{equation}
C(p)=\dfrac{d(\rho\phi S)}{dp}=\phi S\frac{d\rho}{dp}+\rho S\frac{d\phi}{dp}+\rho\phi\frac{dS}{dp},
\end{equation}
in consecuense
\begin{equation}
\dfrac{\partial(\rho\theta)}{\partial t}=C(p)\frac{\partial p}{\partial t}.
\end{equation}
In the classical model, the particles move smoothly in the medium, giving a relationship mean square displacement linear with time, this movement is said to be Markovian because it does not depend on previous positions and is represented by a gaussian diffusion equation . This relationship does not exist in the anomalous diffusion that is altered due to particle entrapment, pore throat closure, fractures, which gives a non-linear relationship, it is said that this is non-Markovian because the movement depends on the previous time or previous position equation \eqref{anomala}.
\begin{equation}\label{normal}
\left\langle x^{2}\right\rangle \propto t
\end{equation}
\begin{equation}\label{anomala}
\left\langle x^{2}\right\rangle \propto t^{\alpha};\:\alpha\neq 1
\end{equation}
in the case that $ 0<\alpha<1 $ it is said to be a subdiffusion, and if $ \alpha>1 $ is a superdiffusion  \cite{Metzler1994}.
The memory formalism can be expressed as a modification of Darcy's Law, through the fractional derivative Caputo see \eqref{caputodef},this equation can explain the diffusion in a case where the effects of past events exist \cite{raghavan2018conceptual}:

%In some sites the structures and  behavior of the fluids are not ideal, emerging the concept of memory, with this concept the behavior of the fluid depends on its space-time trajectory and not of the classic Markovian. The memory concept has been developed by several researchers, including Caputo in
%\cite{Caputo2004}. We have also tried to reflect a fractal structure of the media in the model. In our case, we will use a derivative of the Caputo type for a version of Darcy's law, described by Le Mehaute in \cite{LeMehaute1984} and that Raghavan rewrites in \cite{Raghavan2013} as:
\begin{equation}\label{Darcyfrac}
v(x,t)=-\dfrac{K_{\alpha}}{\mu}\dfrac{\partial^{1-\alpha}}{\partial t^{1-\alpha}}\dfrac{\partial p(x,t)}{\partial x},
\end{equation}

%$ \gamma = \dfrac{1}{d_{f}} $, where $ d_{f} $ as the Hausdorff fractal dimension of the medium, along with the conservation equation in rectangular coordinates:
$ \alpha $ is a parameter, Raghavan and Chen\cite{raghavan2018conceptual} has directly addressed estimates for $ \alpha $: $ 0.56<\alpha<0.91 $ and $ 0.77<\alpha<0.94 $ for fractures and matrix of rock respectively.
Considering cylindrical coordinates and assuming the well axis passes through the origin, the application of the conservation of mass principle to a control volume:
\begin{equation}\label{conmas}
\dfrac{1}{r}\dfrac{\partial}{\partial r}v(r,t)=\phi c\dfrac{\partial}{\partial t}p(r,t)
\end{equation}
where $ \phi $ is the porosity of the medium, and $ r $ is the distance from the well center.
On substituting the left-hand side of equation \ref{conmas} for $ v(r,t) $;   we obtain
the partial differential equation for transient diffusion under subdiffusive flow to be
\begin{equation}
\dfrac{1}{r}\frac{\partial}{\partial r }\left( \dfrac{\tilde{k}}{\mu}\dfrac{1}{r}\dfrac{\partial ^{1-\alpha}}{\partial t^{1-\alpha}}\dfrac{\partial p(r,t)}{\partial r}\right) =\phi c\dfrac{\partial}{\partial t}p(r,t)
\end{equation}
since $ \frac{\tilde{k}}{\mu} $ it does not depend on $ r $, the above equation becomes
\begin{equation}
\dfrac{\partial^{1-\alpha}}{\partial t^{1-\alpha}}\left( \dfrac{1}{r}\frac{\partial}{\partial r }\left( \dfrac{\tilde{k}}{\mu}\dfrac{1}{r}\dfrac{\partial p(r,t)}{\partial r}\right)\right)  =\phi c\dfrac{\partial}{\partial t}p(r,t)
\end{equation}
 then applying the Riemann-Liuoville integral: $ {_t}J^{\alpha-1} $   to both sides, applying the equation \eqref{incaputo} is obtained:
 \begin{equation}\label{Ecprin}
 \dfrac{1}{r}\frac{\partial}{\partial r }\left( \dfrac{\tilde{k}}{\mu}\dfrac{1}{r}\dfrac{\partial p(r,t)}{\partial r}\right)  =\phi c\dfrac{\partial^{\alpha}}{\partial t^{\alpha}}p(r,t)
 \end{equation}

%\begin{equation}
%\dfrac{\partial}{\partial x_{i}} q_{i}(\bar{x};t)=\phi c\dfrac{\partial}{\partial t}p(\bar{x},t),
%\end{equation}
%when combining the two previous equations it is obtained with a system with radial symmetry
%\begin{equation}\label{ecu08}
 %\dfrac{1}{r^{n-1}}\dfrac{\partial}{\partial r} \left[r^{n-1}\lambda(r)\dfrac{\partial p(r,t)}{\partial r}\right] =\phi c\dfrac{\partial^{2-\gamma}}{\partial t^{2-\gamma}}p(r,t),
 %\end{equation}
 %where $ n $ is the Euclidean dimension of the medium, in our case $ n = 2 $.
%Holy \cite{Holy2016} proposes a 1-D linear diffusion model for a horizontal well, based on a generalized Darcy Law, with the use of Caputo-type derivatives of a fractional order in time and space, expressed in the following equation:
%\begin{equation}
%\dfrac{\partial}{\partial x}\left[\dfrac{k_{\alpha \beta}}{\mu}\dfrac{\partial^{1-\alpha}}{\partial t^{1-\alpha}}\left(\dfrac{\partial^{\beta}p(x,t)}{\partial x^{\beta}}\right)\right]=\phi c_{t}\dfrac{\partial p(x,t)}{\partial t}
%\end{equation}
%\begin{equation}
%p(x,0)=p_{i}, 0\leq x\leq L,\;u(L,t)=0,\;0\leq t
%\end{equation}
%Being the conditions of constant rate and constant pressure:
%\begin{equation}
%u(0,t)=\dfrac{q_{f}B}{L_{f}h},\;p(0,t)=p_{f}
%\end{equation}
%The equation is solved with finite differences, and produces the expression:
%\begin{equation}
%\log[p_{i}-p_{f}(t)]=\log(a)-m\log(t)
%\end{equation}
%where $a$ is constant and $m=1-\dfrac{\alpha}{\beta+1}$, 

\subsection{Fractional Calculus}
There are several definitions of fractional derivative: the most widespread is that of Riemann-Liouville, we will only provide the definition of the derivative Caputo because it is the one we will use, a very complete reference in the area can be consulted in the book by Baleanu et al. \cite{Baleanu2011}.
The left-sided fractional integral of Riemann-Liouville of order $ \alpha$ is defined as
\begin{equation}
{_t}J^{\alpha}f(t):=\dfrac{1}{\Gamma(\alpha)}\int_{0}^{t}(t-\tau)^{\alpha-1}f(\tau)d\tau, \:\alpha>0.	
\end{equation}
Where the convention  $ {_t}J^0=I $ (identity operator) and the semigroup property:
\begin{equation}
{_t}J^{\alpha}{_t}J^{\beta}={_t}J^{\beta}{_t}J^{\alpha}={_t}J^{\alpha+\beta}, \:\alpha,\beta\geq0.
\end{equation}
We define the left-sided Caputo fractional derivative  of order $ \alpha> 0 $ as the operator $ {_t}D^{\mu}_{\ast} $ such that $ {_t}D^{\mu}_{\ast}f(t):={_t}J^{m- \mu} {_t}D^{m}f(t) $, with $ D=\frac{d}{dx} $ hence
\begin{equation}\label{caputodef}
{_t}D^{\mu}_{\ast}f(t)\hspace{-.3em}=\hspace{-.3em}\begin{cases}
\hspace{-.3em}\dfrac{1}{\Gamma(m-\mu)} \displaystyle\int_{0}^{t}\dfrac{f^{(m)}(\tau)d\tau}{(t-\tau)^{\mu+1-m}}\hspace{-.3em},\,m-1<\mu<m
\\
\\\hspace{-.3em}\dfrac{d^{m}}{dt^{m}}f(t), \;\mu=m
\end{cases}\hspace{-1em}.
\end{equation}
with $ m=\lfloor \alpha \rfloor $+1 if $ \alpha\notin \mathbb{N} $, $ m=\alpha $ if $ \alpha\in\mathbb{N}. $\\ 
In particular, when $ 0<\alpha<1 $ then
\begin{equation}
{_t}D^{\alpha}_{\ast}={_t}J^{1-\alpha}D
\end{equation}
The above equation implies:
\begin{equation}\label{incaputo}
{_t}J^{\alpha}{_t}D^{\alpha}_{\ast}f(t)={_t}J^{\alpha}{_t}J^{1-\alpha}Df(t)=JDf(t)=f(t)-f(0)
\end{equation}
The fractional derivative Caputo satisfies the property that it is zero when is applied to a constant. Another important property is that you can apply a Laplace transform:
\begin{gather}\label{ecu12}
\mathcal{L}\{{_t}D^{\mu}_{\ast}f(t);s\}=s^{\mu}\tilde{f}(s)-\sum^{m-1}_{k=0}s^{\mu-1-k}f^{(k)}(0^{+}), \hspace{0.2cm} m-1<\mu<k,
\end{gather}
where  

\begin{eqnarray*}
\displaystyle\tilde{f}(s)=\mathcal{L}\{f(t);s\}= \int^{\infty}_{0}e^{-st}f(t)dt, &s\in\mathbb{C},
\end{eqnarray*}

and

\begin{eqnarray*}
f^{(k)}(0^{+}):=\lim\limits_{t \to 0^{+}}f(t).
\end{eqnarray*}

\subsection{Bessel Functions}
The following differential equation of second order
\begin{equation}\label{ecu13}
z^{2}\dfrac{d^{2}y}{dz^{2}}+z\dfrac{dy}{dz}-(z^{2}+\nu^{2})y=0,
\end{equation}
where $ \nu $ is a real constant is called the modified Bessel equation, the solutions to the above equation are called modified Bessel functions which take the following form:
\begin{equation}\label{ecu14}
K_{\nu}(z)=\left(\dfrac{\pi}{2}\right)\dfrac{I_{-\nu}(z)-I_{\nu}(z)}{\mathrm{sin}(\nu\pi)},
\end{equation}
where $ I_{\nu}(z) $ are the modified Bessel functions of the first type, it is noted that $ I_{\nu} $ and $ I_{-\nu}$ - form a set of solutions for the equation \eqref{ecu13} and the equation \eqref{ecu14} it is known as the modified Bessel function of the second type.
Some properties of the modified Bessel function of the second type are:
\begin{equation}\label{ecu15}
\dfrac{d}{dz}K_{\nu}(\alpha z)=-\alpha K_{\nu-1}(\alpha z)-\dfrac{\nu}{z}K_{\nu}(\alpha z),
\end{equation}
\begin{equation}\label{ecu16}
\dfrac{d}{dz}K_{\nu}(\alpha z)=-\alpha K_{\nu+1}(\alpha z)+\dfrac{\nu}{z}K_{\nu}(\alpha z).
\end{equation}
\subsection{Flow equation  (fractional time derivative)}
Let's assume the equation \eqref{Ecprin} that represents the fluid, where the medium is a whole, so we have:
\begin{equation}\label{ecu17}
\phi c_{\alpha}\dfrac{\partial^{\alpha}p}{\partial t^{\alpha}}=\dfrac{k}{\mu}\dfrac{1}{r}\dfrac{\partial}{\partial r}\left(r\dfrac{\partial p}{\partial r}\right),
\end{equation}
where the derivative expression on the left denotes the fractional derivative Caputo of order $ \alpha\ \in \mathbb {R} $, with dimensionless variables, the equation \eqref{ecu17} is
\begin{equation}\label{ecu18}
\phi  c_{D\alpha}\dfrac{\partial^{\alpha}p_{D}}{\partial t^{\alpha}_{D}}=\kappa_{D}\dfrac{1}{r}\dfrac{\partial}{\partial r_{D}}\left(r_{D}\dfrac{\partial p_{D}}{\partial r_{D}}\right),
\end{equation}
where
\begin{equation}
p_{D}=\dfrac{2\pi hk(p_{i}-p)}{Q_{0}B_{0}\mu},\hspace{0.2cm} \:t_{D}=t\dfrac{k}{\phi cr^{2}_{w}\mu},\hspace{0.2cm} \:r_{D}=\dfrac{r}{r_{w}},
\end{equation}
where $ \phi $ represents the porosity (dimensionless), $ c $ represents the compressibility of the medium in units of $ Pa^{-1 } $, $ k $ represents the permeability of the medium with units of $ m^{2} $, $ p $ represents the fluid pressure in the middle with units of $ Pa $, $ \mu $ is the viscosity of the fluid with units of $ Pa\cdot s $,  $t$ represents the time in units of $s$, $r$ represents the distance of the well in units of $ m $, $ r_ {w} $ is a reference parameter: well radius with units of $ m $, $ h $ is the well thickness measured in $ m $, $ p_{i} $ is the initial reservoir pressure, the value of $ Q_{0} $ is the flow rate with units of $ m^{3}  s^{-1}$ and $ B_{0} $ is the fluid factor (dimensionless).
\subsection{ Laplace Transform}
The Laplace transform applied to the equation \eqref{ecu18} gives the following result using the equation \eqref{ecu12}
\begin{equation}\label{ecu20}
u^{\alpha}\bar{p}_{D}=\dfrac{1}{r_{D}}\dfrac{\partial}{\partial r_{D}}\left(r_{D}\dfrac{\partial\bar{p}_{D}}{\partial r_{D}}\right),\ u>0,
\end{equation}
where $ \bar{p}_D(t_{0})=0 $, because $ p=p_{i} $, in $ t=t_{0} $.
\subsubsection{ Bessel Functions}
The spatial derivatives to be developed in the equation \eqref{ecu20} present the following form:
\begin{equation}
r^{2}_{D}\dfrac{\partial^{2}\bar{p}_{D}}{\partial r^{2}_{D}}+r_{D}\dfrac{\partial\bar{p}_{D}}{\partial r_{D}}-r^{2}_{D}u^{\alpha}\bar{p}_{D}=0,
\end{equation}
which is a Bessel equation, so the solution is:
 \begin{equation}\label{ecu22}
\bar{p}_{D}=AK_{0}(\beta r_{D}).
\end{equation}
By substituting the equation \eqref{ecu22} into the equation \eqref{ecu20} and considering  the equations \eqref{ecu15} and \eqref{ecu16} to find the value of $ \beta $ we have :
\begin{equation}\label{ecu23}
\beta=\pm\sqrt{u^{\alpha}}, \:u>0.
\end{equation}
The equation \eqref{ecu22} when considering the value of $ \beta$, equation \eqref{ecu23} is
\begin{equation}\label{ecu24}
\bar{p}_{D}=AK_{0}(r_{D}\sqrt{u^{\alpha}}).
\end{equation}
In the equation \eqref{ecu24} $ \beta=-\sqrt{u^{\alpha}} $ is discarded because the modified Bessel function of second type is not defined for negative values
\subsubsection{Border conditions}
To find the solution to the equation \eqref{ecu18}, the following boundary condition is considered:
\begin{equation}\label{ecu25}
r_{D}\dfrac{\partial\bar{p}_{D}}{\partial r_{D}}\Big|_{r_{D}=1}=-\dfrac{1}{u}.
\end{equation}
Substituting the equation \eqref{ecu24} in \eqref{ecu25} generates the following:
 \begin{equation}
A=\dfrac{1}{u}[\sqrt{u^{\alpha}}K_{1}(\sqrt{u^{\alpha}})]^{-1},
\end{equation}
\begin{equation}
\bar{p}_{D}=\dfrac{1}{u}[\sqrt{u^{\alpha}}K_{1}(\sqrt{u^{\alpha}})]^{-1}K_{0}(r_{D}\sqrt{u^{\alpha}}).
\end{equation}
Therefore, the value of the pressure at the boundary of the well ($ r_{D}=1 $) is in the space of Laplace:
 \begin{equation}
\bar{p}_{D}|_{r_{D}=1}=\dfrac{1}{u}[\sqrt{u^{\alpha}}K_{1}(\sqrt{u^{\alpha}})]^{-1}K_{0}(\sqrt{u^{\alpha}}).
\end{equation}

\section{Flow equation with triple porosity and triple permeability with fractional time derivative}
From the classical transfer equations, B. Martínez in \cite{martinez2017applications2} proposes a system of coupled flow equations with triple porosity and triple permeability, which have the following form: \\
We consider the following notation
\begin{equation*}
c_{s_1}:=\dfrac{1}{\phi_{s_1}}\dfrac{\partial\phi_{s_1}}{\partial p_{s_1}},\text{ with }\,s_1=m,f,v
\end{equation*}

\begin{equation*}
\Delta_{s_1 s_2}(p):=p_{s_1}-p_{s_2}	
\end{equation*}

\begin{equation}\label{ecu29}
\phi_{m}c_{m}\dfrac{\partial p_{m}}{\partial t}=\dfrac{k_{m}}{\mu}\dfrac{1}{r}\dfrac{\partial}{\partial r}\left(r\dfrac{\partial p_{m}}{\partial r}\right)+a_{mf}\Delta_{fm}(p)+a_{mv}\Delta_{vm}(p)
\end{equation}
\begin{equation}\label{ecu30}
\phi_{f}c_{f}\dfrac{\partial p_{f}}{\partial t}=\dfrac{k_{f}}{\mu}\dfrac{1}{r}\dfrac{\partial}{\partial r}\left(r\dfrac{\partial p_{f}}{\partial r}\right)-a_{mf}\Delta_{fm}(p)+a_{fv}\Delta_{vf}(p)
\end{equation}
\begin{equation}\label{ecu31}
\phi_{v}c_{v}\dfrac{\partial p_{v}}{\partial t}=\dfrac{k_{v}}{\mu}\dfrac{1}{r}\dfrac{\partial}{\partial r}\left(r\dfrac{\partial p_{v}}{\partial r}\right)-a_{mv}\Delta_{vm}(p)-a_{fv}\Delta_{vf}(p)
\end{equation}

\begin{comment}
\begin{itemize}
\item[•] If $a_{s_1s_2}=a_{s_2s_1} \ \forall s_1,s_2\in \set{m,f,v}$

\begin{equation}
\phi_{m}c_{m}\dfrac{\partial p_{m}}{\partial t}=\dfrac{k_{m}}{\mu}\dfrac{1}{r}\dfrac{\partial}{\partial r}\left(r\dfrac{\partial p_{m}}{\partial r}\right)-a_{mf}\Delta_{mf}(p)-a_{mv}\Delta_{mv}(p)
\end{equation}
\begin{equation}
\phi_{f}c_{f}\dfrac{\partial p_{f}}{\partial t}=\dfrac{k_{f}}{\mu}\dfrac{1}{r}\dfrac{\partial}{\partial r}\left(r\dfrac{\partial p_{f}}{\partial r}\right)-a_{fm}\Delta_{fm}(p)-a_{fv}\Delta_{fv}(p)
\end{equation}
\begin{equation}
\phi_{v}c_{v}\dfrac{\partial p_{v}}{\partial t}=\dfrac{k_{v}}{\mu}\dfrac{1}{r}\dfrac{\partial}{\partial r}\left(r\dfrac{\partial p_{v}}{\partial r}\right)-a_{vm}\Delta_{vm}(p)-a_{vf}\Delta_{vf}(p)
\end{equation}

then if $s_1\neq s_2 \neq s_3$

\begin{eqnarray}
\phi_{s_1}c_{s_1}\dfrac{\partial p_{s_1}}{\partial t}=\dfrac{k_{s_1}}{\mu}\dfrac{1}{r}\dfrac{\partial}{\partial r}\left(r\dfrac{\partial p_{s_1}}{\partial r}\right)-a_{s_1s_2}\Delta_{s_1s_2}(p)-a_{s_1s_3}\Delta_{s_1s_3}(p), &\forall s_1,s_2,s_3\in \set{m,f,v}
\end{eqnarray}

\end{itemize}

\end{comment}

where $ \phi_{m}, \phi_{f}, \phi_{v} $ represent the porosities of the soil matrix, the fractured medium and the vugular medium respectively in units of $ m^{3} / m^{ 3} $; $ c_{m}, c_{f}, c_{v} $ represent the compressibility in each porous medium in units of $ Pa^{-1} $; $ k_{m}, k_{f}, k_{v} $ represent the permeability of each porous medium with units of $ m^{2} $; $ p_{m}, p_{f}, p_{v} $ represent the fluid pressure in each medium porous with units of $ Pa $; $ \mu $ is the viscosity of the fluid with units of $ Pa \cdot s $; $ a_{mf}, a_{mv}, a_{fv} $ are the transfer terms in the matrix-fracture, matrix-void, and fracture-vortex interfaces respectively with units of $ Pa^{-1}s^{-1} $, $ t $ represents the time in units of $ s $ and $ r $ represents the distance to the well in units of $ m $.
\subsection{Adimensionalization of the flow equations}
In order to handle the equations \eqref{ecu29}, \eqref{ecu30} and \eqref{ecu31} in an easier way, the dimensionlessness of variables is applied. The dimensionlessness is a technique commonly used to make the parameters or variables in an equation have no units, to rank the possible values of a variable or a constant in order that its value is known and thus more manipulable. 
The system of equations \eqref{ecu29}, \eqref{ecu30} and \eqref{ecu31} takes, after applying the dimensionlessness, the following form
\begin{equation}\label{ecu32}
\omega_m \dfrac{\partial p_{Dm}}{\partial t_{D}}=\kappa_m \dfrac{1}{r_{D}}\dfrac{\partial}{\partial r_{D}}\left(r_{D}\dfrac{\partial p_{Dm}}{\partial r_{D}}\right)
+\lambda_{mf}\Delta_{fm}(p_{D})+\lambda_{mv}\Delta_{vm} (p_{D}),
\end{equation}

\begin{equation}\label{ecu33}
\omega_{f}\dfrac{\partial p_{Df}}{\partial t_{D}}=\kappa_{f}\dfrac{1}{r_{d}}\dfrac{\partial}{\partial r_{D}}\left(r_{D}\dfrac{\partial p_{Df}}{\partial r_{D}}\right)
-\lambda_{mf}\Delta_{fm} (p_{D})+\lambda_{fv}\Delta_{vf} (p_{D}),
\end{equation}

\begin{equation}\label{ecu34}
\omega_{v}\dfrac{\partial p_{Dv}}{\partial t_{D}}=\kappa_{v}\dfrac{1}{r_{d}}\dfrac{\partial}{\partial r_{D}}\left(r_{D}\dfrac{\partial p_{Dv}}{\partial r_{D}}\right)
-\lambda_{mv}\Delta_{vm} (p_{D})-\lambda_{fv}\Delta_{vf} (p_{D}),
\end{equation}

with

\begin{align}
\omega_m=&1-\omega_f-\omega_v,\\
\kappa_m=&1-\kappa_f-\kappa_v,
\end{align}

where

	\begin{align}\label{ecu35}
\omega_{s_1}=\dfrac{\phi_{s_1}c_{s_1}}{\phi_{m}c_{m}+\phi_{f}c_{f}+\phi_{v}c_{v}}, \hspace{0.2cm} s_1=f,v,
	\end{align}

	\begin{align}\label{ecu36}
	r_{D}=\dfrac{r}{r_{w}}, \hspace{0.2cm} \kappa_{s_1}=\dfrac{k_{s_1}}{k_{m}+k_{f}+k_{v}}, \hspace{0.2cm} s_1=f,v,
	\end{align}

	\begin{align}\label{ecu37}
\lambda_{s_1 s_2}=\dfrac{a_{s_1 s_2}\mu r^{2}_{w}}{k_{m}+k_{f}+k_{v}}, \hspace{0.2cm} s_1 s_2=mf,mv,fv,
	\end{align}

\begin{subequations}\label{ecu38}
	\begin{align}
	&p_{Dj}=\dfrac{2\pi h(k_{m}+k_{f}+k_{v})(p_{i}-p_{j})}{Q_{0}B_{0}\mu},\\&t_{D}=\dfrac{t(k_{m}+k_{f}+k_{v})}{\mu r^{2}_{w}(\phi_{m}c_{m}+\phi_{f}c_{f}+\phi_{v}c_{v})}.
	\end{align}
\end{subequations}
The equations \eqref{ecu35} - \eqref{ecu38} represent the dimensionless variables, it can be verified that these variables do not have units; in the equation \eqref{ecu37} the value of $ r_{w} $ is a reference parameter, in this case the radius of the well, in order that the variable $ r_{D} $ has the minimum value equal to $ 1 $, with units of $ m $. In the equation \eqref{ecu38} the value of $ h $ represents the thickness of the oil field with units of $ m $; $ p_{j} $ is the pressures in the different porous media, where $ j = m, f, v, p_{i}; $ is the initial pressure in the field; the value of $ Q_{0} $ is the flow rate with units of $ m^{3}s^{-1}  $ and $ B_{0} $ is the fluid formation factor (dimensionless).

\subsection{The system with fractional derivative}

From the system of equations with dimensionless variables \eqref{ecu32} - \eqref{ecu34}, using the equation of flow with fractional time derivative \eqref{ecu18}, we express a system with a fractional time derivative:
\begin{align}\label{ecu39}
\omega_m\dfrac{\partial^{\beta_m} p_{Dm}}{\partial t^{\beta_m}_{D}}=\kappa_m \dfrac{1}{r_{D}}\dfrac{\partial}{\partial r_{D}}\left(r_{D}\dfrac{\partial p_{Dm}}{\partial r_{D}}\right)+\lambda_{mf}\Delta_{fm} (p_{D})+\lambda_{mv}\Delta_{vm} (p_{D}),\\\label{ecu40}
\omega_{f}\dfrac{\partial^{\beta_f} p_{Df}}{\partial t^{\beta_f}_{D}}=\kappa_{f}\dfrac{1}{r_{d}}\dfrac{\partial}{\partial r_{D}}\left(r_{D}\dfrac{\partial p_{Df}}{\partial r_{D}}\right)-\lambda_{mf}\Delta_{fm} (p_{D})+\lambda_{fv}\Delta_{vf} (p_{D}),\\\label{ecu41}
\omega_{v}\dfrac{\partial^{\beta_v} p_{Dv}}{\partial t^{\beta_v}_{D}}=\kappa_{v}\dfrac{1}{r_{d}}\dfrac{\partial}{\partial r_{D}}\left(r_{D}\dfrac{\partial p_{Dv}}{\partial r_{D}}\right)-\lambda_{mv}\Delta_{vm} (p_{D})-\lambda_{fv}\Delta_{vf}(p_{D}),
\end{align}
where the variables shown in the equations \eqref{ecu39} - \eqref{ecu41} have the same meaning as the equations \eqref {ecu35} - \eqref{ecu37}. By means of the Laplace transform and with the use of the equation \eqref{ecu12} the following system is reached:
\begin{align}\label{ecu42}
	\omega_m u^{\beta_m}\bar{p}_{Dm}=\kappa_m \dfrac{1}{r_{D}}\dfrac{\partial}{\partial r_{_D}}\left(r_{D}\dfrac{\partial \bar{p}_{Dm}}{\partial r_{d}}\right)+\lambda_{mf}\Delta_{fm} (\bar{p}_{D})+\lambda_{mv}\Delta_{vm} (\bar{p}_{D}),\\\label{ecu43}
	\omega_{f}u^{\beta_f}\bar{p}_{Df}=\kappa_{f}\dfrac{1}{r_{D}}\dfrac{\partial}{\partial r_{_D}}\left(r_{D}\dfrac{\partial \bar{p}_{Df}}{\partial r_{d}}\right)-\lambda_{mf}\Delta_{fm} (\bar{p}_{D})+\lambda_{fv}\Delta_{vf}(\bar{p}_{D}),\\\label{ecu44}
	\omega_{v}u^{\beta_v}\bar{p}_{Dv}=\kappa_{v}\dfrac{1}{r_{D}}\dfrac{\partial}{\partial r_{_D}}\left(r_{D}\dfrac{\partial \bar{p}_{Dv}}{\partial r_{d}}\right)-\lambda_{mv}\Delta_{vm} (\bar{p}_{D})-\lambda_{fv}\Delta_{vf}(\bar{p}_{D}),
\end{align}
where the variables shown in the equations \eqref{ecu39} - \eqref{ecu41} have the same meaning as those where $\bar{p}_{Dm},\bar{p}_{Df}$ and $\bar{p}_{Dv}$ represent the Laplace transforms of the variables $p_{Dm}, p_{Df}$ and $p_{Dv}$. After developing  equations \eqref{ecu42} - \eqref{ecu44} it is easy to see that they comply with the form of a Bessel equation and therefore their solutions, as in the case $\beta_m=\beta_f=\beta_v=1$, they are

\begin{align}\label{ecu45}
&\bar{p}_{Dm}=AK_{0}(\alpha r_{D}),\\\label{ecu46}
&\bar{p}_{Df}=BK_{0}(\alpha r_{D}),\\\label{ecu47}
&\bar{p}_{Dv}=CK_{0}(\alpha r_{D}),
\end{align}
In order to simplify the successive equations, the following terms are defined:

\begin{subequations}\label{ecu48}
	\begin{align}
	m_{1}(u)&=u^{\beta_m}\omega_m+\lambda_{mf}+\lambda_{mv},\\
	m_{2}&=\lambda_{mf},\\
	 m_{3}&=\lambda_{mv},
	\end{align}
\end{subequations}

\begin{subequations}\label{ecu49}
	\begin{align}
	m_{4}(u)&=u^{\beta_f}\omega_{f}+\lambda_{mf}+\lambda_{fv},\\
	m_{5}&=\lambda_{fv},\\
	m_{6}(u)&=u^{\beta_v}\omega_{v}+\lambda_{mv}+\lambda_{fv}.
	\end{align}
\end{subequations}

As a result of replacing the equations \eqref{ecu45} - \eqref{ecu47} in the system shown in \eqref{ecu42}-\eqref{ecu44} and making use of the definitions shown by \eqref{ecu48}-\eqref{ecu49}, we have the following:

\begin{align}\label{ecu50}
&K_{0}(\alpha r_{D})\{A[\kappa_m\alpha^{2}-m_{1}]+Bm_{2}+Cm_{3}\}=0,\\\label{ecu51}
&K_{0}(\alpha r_{D})\{Am_{2}+B[\kappa_{f}\alpha^{2}-m_{4}]+Cm_{5}\}=0,\\\label{ecu52}
&K_{0}(\alpha r_{D})\{Am_{3}+Bm_{5}+C[\kappa_{v}\alpha^{2}-m_{6}]\}=0.
\end{align}

Since the modified second-species Bessel functions have an asymptotic behavior, that is, they never take the value of zero, then the system shown in the equations \eqref{ecu50}-\eqref{ecu52} can be expressed as follows:

\begin{equation}\label{ecu53}
\begin{bmatrix}
\kappa_m\alpha^2-m_{1}&m_2&m_3\\
m_2&\kappa_f\alpha^2-m_4&m_5\\
m_3&m_5&\kappa_v\alpha^2-m_6
\end{bmatrix}\hspace{-3pt}
\begin{bmatrix}
A\\
B\\
C
\end{bmatrix}=\begin{bmatrix}
0\\
0\\
0
\end{bmatrix}\normalsize
\end{equation}

The equation \eqref{ecu53} is used to find the values of $A,B$ and $C$. Keeping this in mind, we  note two principal cases: the determinant of the matrix $ 3\times3 $ or it is different from zero.
The first case gives us the trivial solution $A=B=C=0$. Linear algebra is known to have a solution to the equation if and only if the determinant equals zero. 
The second case, where the determinant is equal to zero, we obtain the equation of degree six that follows:

\begin{align}
\begin{array}{c}
\kappa_m\kappa_{f}\kappa_{v}\alpha^{6}-[\kappa_m(\kappa_{f}m_{6}+\kappa_{v}m_{4})+\kappa_{f}\kappa_{v}m_{1}]\alpha^{4} \vspace{0.1cm}\\
+[\kappa_mm_{4}m_{6}-\kappa_mm^{2}_{5}+(\kappa_{f}m_{6}+\kappa_{v}m_{4})m_{1}-\kappa_{v}m_{2}^{2}-\kappa_{f}m^{3}_{2}]\alpha^{2} \vspace{0.1cm}\\
-m_{1}m_{4}m_{6}+m_{1}m^{2}_{5}+m^{2}_{2}m_{6}+2m_{2}m_{3}m_{5}+m^{2}_{3}m_{4}=0.
\end{array}
\end{align}

In the previous equation the powers of $\alpha$ are even, therefore it can be solved as an equation of degree three. This equation has three real roots. The general solutions of  equations \eqref{ecu42}-\eqref{ecu44}, when incorporating the three real roots, are:
\begin{align}\label{ecu55}
\bar{p}_{Dm}&=A_{1}D_{1}K_{0}(\alpha_{1}r_{D})+A_{2}D_{2}K_{0}(\alpha_{2}r_{D})+A_{3}D_{3}K_{0}(\alpha_{3}r_{D}),\\\label{ecu56}
\bar{p}_{Df}&=B_{1}D_{1}K_{0}(\alpha_{1}r_{D})+B_{2}D_{2}K_{0}(\alpha_{2}r_{D})+B_{3}D_{3}K_{0}(\alpha_{3}r_{D}),\\\label{ecu57}
\bar{p}_{Dv}&=D_{1}K_{0}(\alpha_{1}r_{D})+D_{2}K_{0}(\alpha_{2}r_{D})+D_{3}K_{0}(\alpha_{3}r_{D}),
\end{align}
where the terms $A_{i}, B_{i},\mathrm{with}\: i=1,2,3$, have the form
\begin{align}
&A_{i}=\dfrac{m_{3}(\kappa_{f}\alpha^{2}_{i}-m_{4})-m_{2}m_{5}}{m^{5}_{2}-[\kappa_m\alpha_{i}^{2}-m_{1}][\kappa_{f}\alpha^{2}_{i}-m_{4}]},\\
&B_{i}=\dfrac{-m_{3}-A_{i}[\kappa_m\alpha^{2}_{i}-m_{1}]}{m_{2}},
\end{align}
where the terms $ D_{1}, D_{2} $ and $ D_{3} $ are obtained from the boundary conditions and $\alpha_{1}, \alpha_{2}, \alpha_{3} $ are the positive roots of  $ \alpha^{2}_{1},\alpha^{2}_{2},\alpha^{2}_{3} $.
The values of $D_{1}, D_{2}, D_{3} $ are obtained from the boundary conditions and they are equal to

\begin{subequations}
	\begin{align}
	\scriptsize
	\begin{array}{c}
	D_1=\dfrac{1}{u}\left\{ \alpha_1E_1K_1(\alpha_1)+\alpha_{3}E_3K_1(\alpha_3)\dfrac{(B_1-1)K_0(\alpha_1)}{(1-B_3)K_0(\alpha_3)}+\dfrac{\left[(1-A_1)K_0(\alpha_{1})+(1-A_3)\dfrac{(B_1-1)}{(1-B_3)}K_0(\alpha_1)\right]}{\left[(A_2-1)K_0(\alpha_{2})+(A_3-1)\dfrac{(B_2-1)}{(1-B_3)}K_0(\alpha_2)\right]}\alpha_{2}E_2K_1(\alpha_{2})+\alpha_{3}E_3K_1(\alpha_3)\dfrac{(B_2-1)K_0(\alpha_2)}{(1-B_3)K_0(\alpha_3)}\right\}^{-1},
	\end{array}
	\end{align}
	
	\begin{align}
	\scriptsize
	\begin{array}{c}
	D_2=\dfrac{1}{u}\left\{\dfrac{\left[(A_2-1)K_0(\alpha_{2})+(A_3-1)\dfrac{(B_2-1)}{(1-B_3)}K_0(\alpha_{2})\right]}{\left[(1-A_1)K_0(\alpha_1)+(1-A_3)\dfrac{(B_1-1)}{(1-B_3)}K_0(\alpha_{1})\right]}+\left[\alpha_{1}E_1K_1(\alpha_{1})+\alpha_{3}E_3K_1(\alpha_{3})\dfrac{(B_1-1)K_0(\alpha_{1})}{(1-B_3)K_0(\alpha_{3})}\right]+\left[\alpha_{2}E_2K_1(\alpha_{2})+\alpha_{3}E_3K_1(\alpha_{3})\dfrac{(B_2-1)K_0(\alpha_{2})}{(1-B_3)K_0(\alpha_{3})}\right]\right\}^{-1},
	\end{array}
	\end{align}
		
	\begin{align}
	\scriptsize
	\begin{array}{c}
	D_3=\dfrac{1}{u}\left\{\alpha_{1}E_1K_1(\alpha_{1})\dfrac{(1-B_3)K_(\alpha_{3})}{(B_1-1)K_0(\alpha_{1})}+\alpha_{3}E_3K_1(\alpha_{3})+\left[\dfrac{(1-A_1)(1-B_3)+(1-A_3)(B_1-1)}{(A_2-1)(1-B_3)+(A_3-1)(B_2-1)}\right]\cdot \left[\dfrac{\alpha_{2}E_2(1-B_3)K_1(\alpha_{2})K_0(\alpha_{3})+\alpha_{3}E_3(B_2-1)K_1(\alpha_{3})K_0(\alpha_{2})}{(B_1-1)K_0(\alpha_{2})}\right]\right\}^{-1} \vspace{0.1cm} \\
	+\dfrac{1}{u}\left\{\alpha_{2}E_2K_1(\alpha_{2})\dfrac{(1-B_3)K_0(\alpha_{3})}{(B_2-1)K_0(\alpha_{2})}+\alpha_{3}E_3K_1(\alpha_{3})+\left[\dfrac{(1-A_2)(1-B_3)+(1-A_3)(B_2-1)}{(A_1-1)(1-B_3)+(A_3-1)(B_1-1)}\right]\cdot
	 \left[\dfrac{\alpha_{1}E_1(1-B_3)K_1(\alpha_{1})K_0(\alpha_{3})+\alpha_{3}E_3(B_1-1)K_1(\alpha_{3})K_0(\alpha_{1})}{(B_2-1)K_0(\alpha_{1})}\right]\right\}^{-1},
	 \end{array}
	\end{align}
\end{subequations}
where
\begin{align*}
E_i=\left[\kappa_mA_i+\kappa_{f}B_i+\kappa_{v}\right],\hspace{0.2cm} i=1,2,3.
\end{align*}
\normalsize
The following equations are obtained at as the result of substituting the equations \eqref{ecu55}-\eqref{ecu57} in boundary conditions.
\begin{align}\label{ecu65}
\alpha_{1}K_{1}(\alpha_{1})D_{1}[\kappa_mA_{1}+\kappa_{f}B_{1}+\kappa_{v}]+\alpha_{2}K_{1}(\alpha_{2})D_{2}[\kappa_mA_{2}+\kappa_{f}B_{2}+\kappa_{v}]\nonumber\\+\alpha_{3}K_{1}(\alpha_{3})D_{3}[\kappa_mA_{3}+\kappa_{f}B_{3}+\kappa_{v}]&=\dfrac{1}{u},\\
(A_{1}-1)D_{1}K_{0}(\alpha_{1})+(A_{2}-1)D_{2}K_{0}(\alpha_{2})+(A_{3}-1)D_{3}K_{0}(\alpha_{3})&=0 \label{ecu66}\\
(B_{1}-1)D_{1}K_{0}(\alpha_{1})+(B_{2}-1)D_{2}K_{0}(\alpha_{2})+(B_{3}-1)D_{3}K_{0}(\alpha_{3})&=0.\label{ecu67}
\end{align}
In order to simplify the equations \eqref{ecu65}-\eqref{ecu67} the following terms are defined:
\begin{align}
&P_{i}=\alpha_{i}K_{1}(\alpha_{i})[\kappa_mA_{i}+\kappa_{f}B_{i}+\kappa_{v}],\\
&Q_{i}=(A_{i}-1)K_{0}(\alpha_{i}),\;R_{i}=(B_{i}-1)K_{0}(\alpha_{i}),
\end{align}
where $i=1,2,3$. The matrix equation associated with the system of equations \eqref{ecu65}-\eqref{ecu67} has the form
\begin{equation}
\begin{bmatrix}
P_1&P_2&P_3\\
Q_1&Q_2&Q_3\\
R_1&R_2&R_3
\end{bmatrix}\hspace{-5pt}
\begin{bmatrix}
D_1\\
D_2\\
D_3
\end{bmatrix}=\begin{bmatrix}
1/u\\
0\\
0
\end{bmatrix}.
\end{equation}
It defines
\begin{align}
m&=Q_{1}R_{2}P_{3}-Q_{1}P_{2}R_{3}-R_{1}Q_{2}P_{3}-R_{2}P_{1}Q_{3}+P_{2}R_{1}Q_{3}+P_{1}Q_{2}R_{3}.
\end{align}
The solution of the matrix equation is
\begin{equation}
\begin{bmatrix}
D_1\\D_2\\D_3
\end{bmatrix}=
\dfrac{1}{m}\begin{bmatrix}
Q_1\\Q_2\\Q_3
\end{bmatrix}\times
\begin{bmatrix}
R_1\\R_2\\R_3
\end{bmatrix}
\end{equation}
Now that we have obtained all the terms that define the solution of the system of equations \eqref{ecu42}-\eqref{ecu44} have been obtained, it is necessary to show the value of the pressure at the boundary of the well, this value is:

\begin{align}\label{eq:001} 
\bar{p}_{w}\big{|}_{r_{d}=1}= \sum_{i=1}^3 D_iK_0(\alpha_i) =\sum_{i=1}^3A_i D_iK_0(\alpha_i)=\sum_{i=1}^3B_i D_iK_0(\alpha_i).
\end{align}

\begin{figure}[!ht]
     \centering
     \begin{subfigure}[b]{0.49\textwidth}
         \centering
         \includegraphics[width=\textwidth]{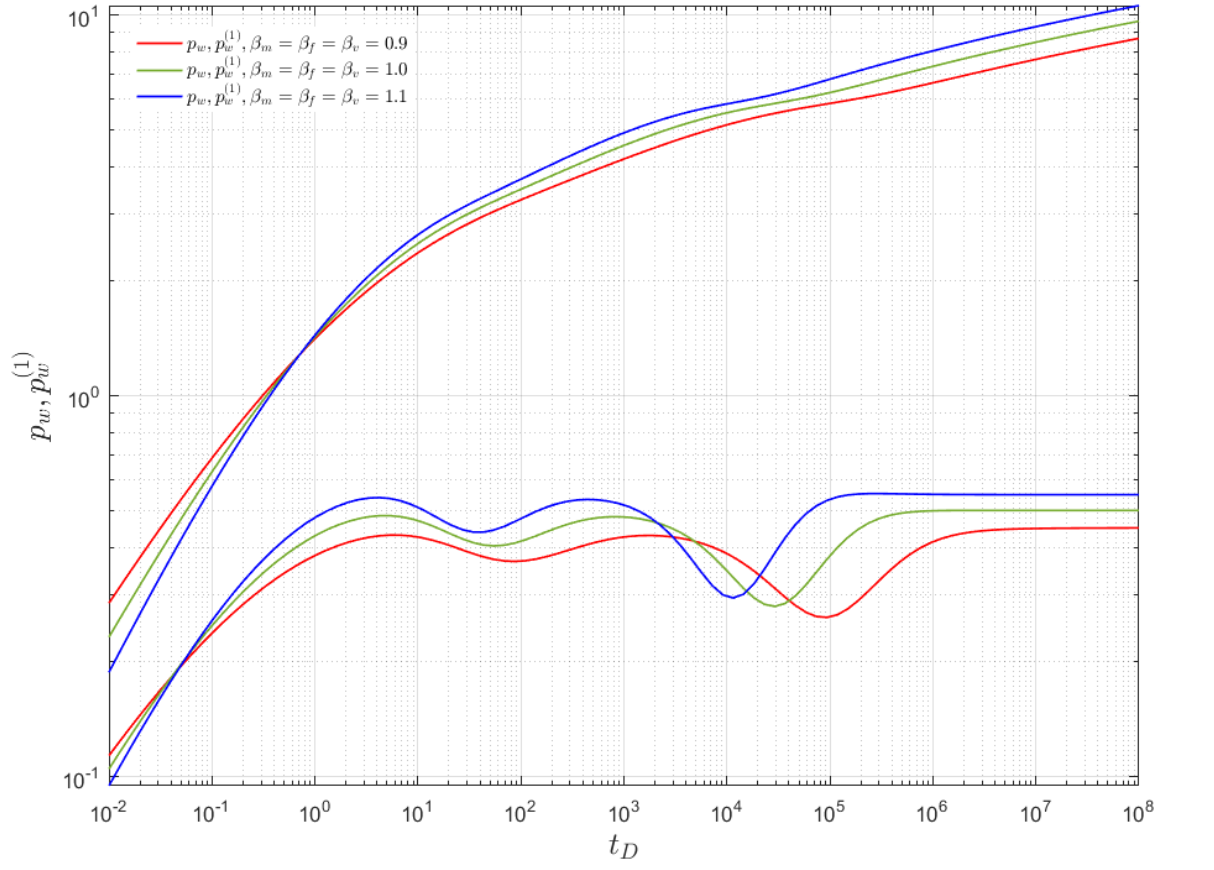}
     \end{subfigure}
     \begin{subfigure}[b]{0.49\textwidth}
         \centering
         \includegraphics[width=\textwidth]{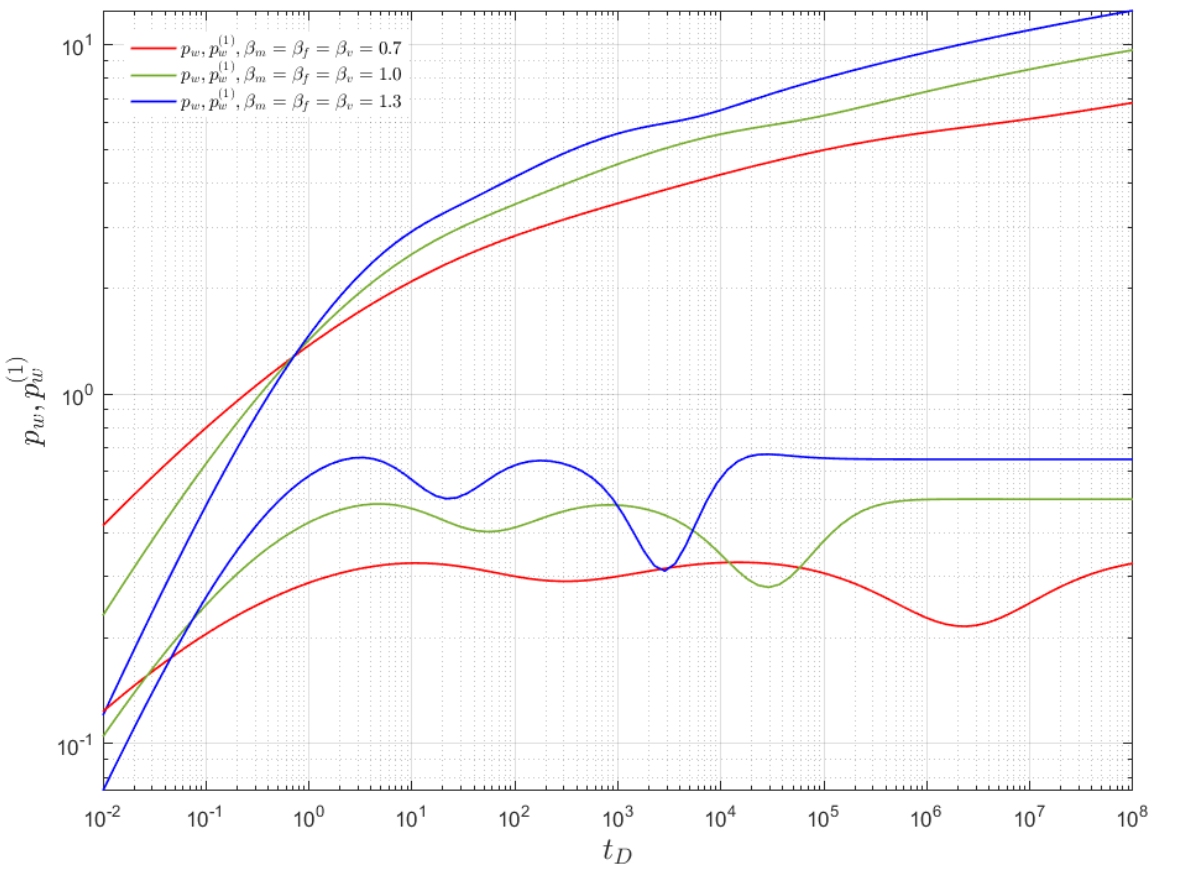}
     \end{subfigure}
     \begin{subfigure}[b]{0.49\textwidth}
         \centering
         \includegraphics[width=\textwidth]{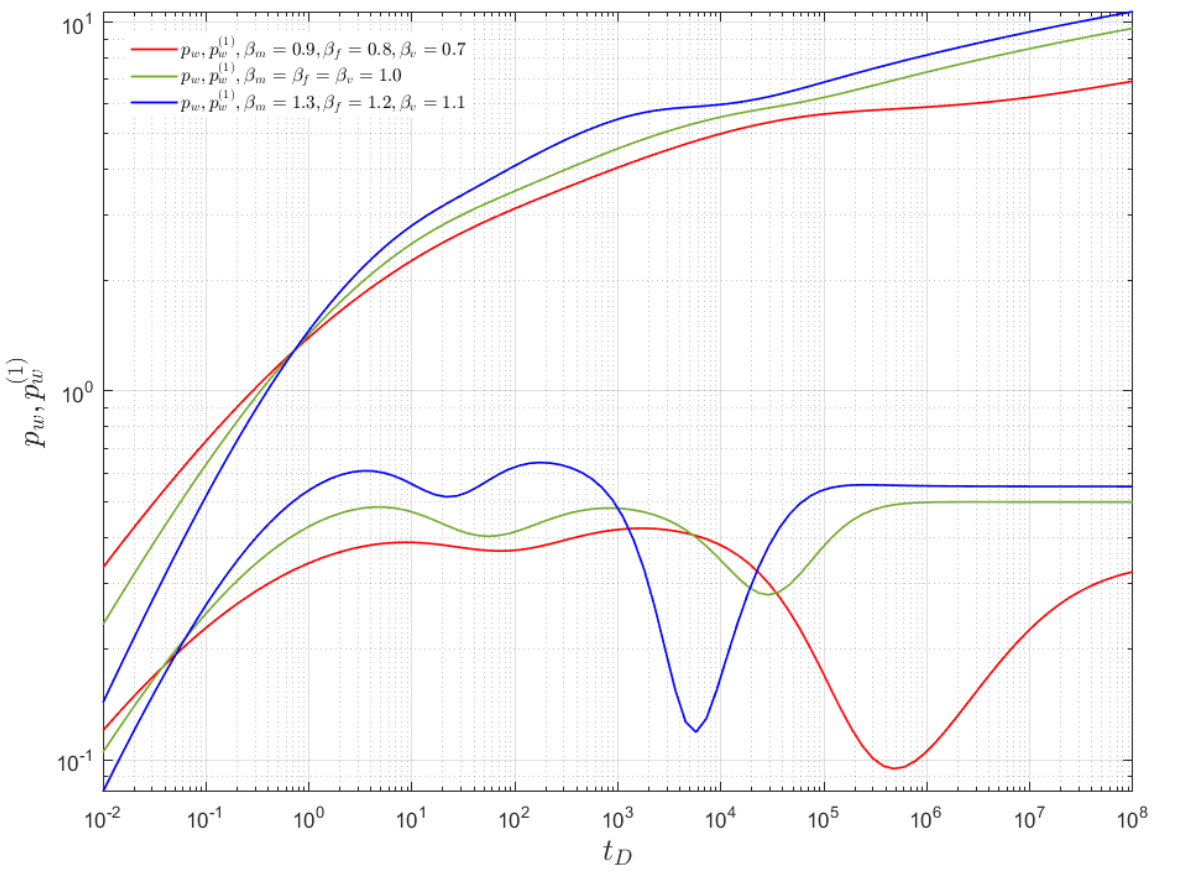}
     \end{subfigure}
          \begin{subfigure}[b]{0.49\textwidth}
         \centering
         \includegraphics[width=\textwidth]{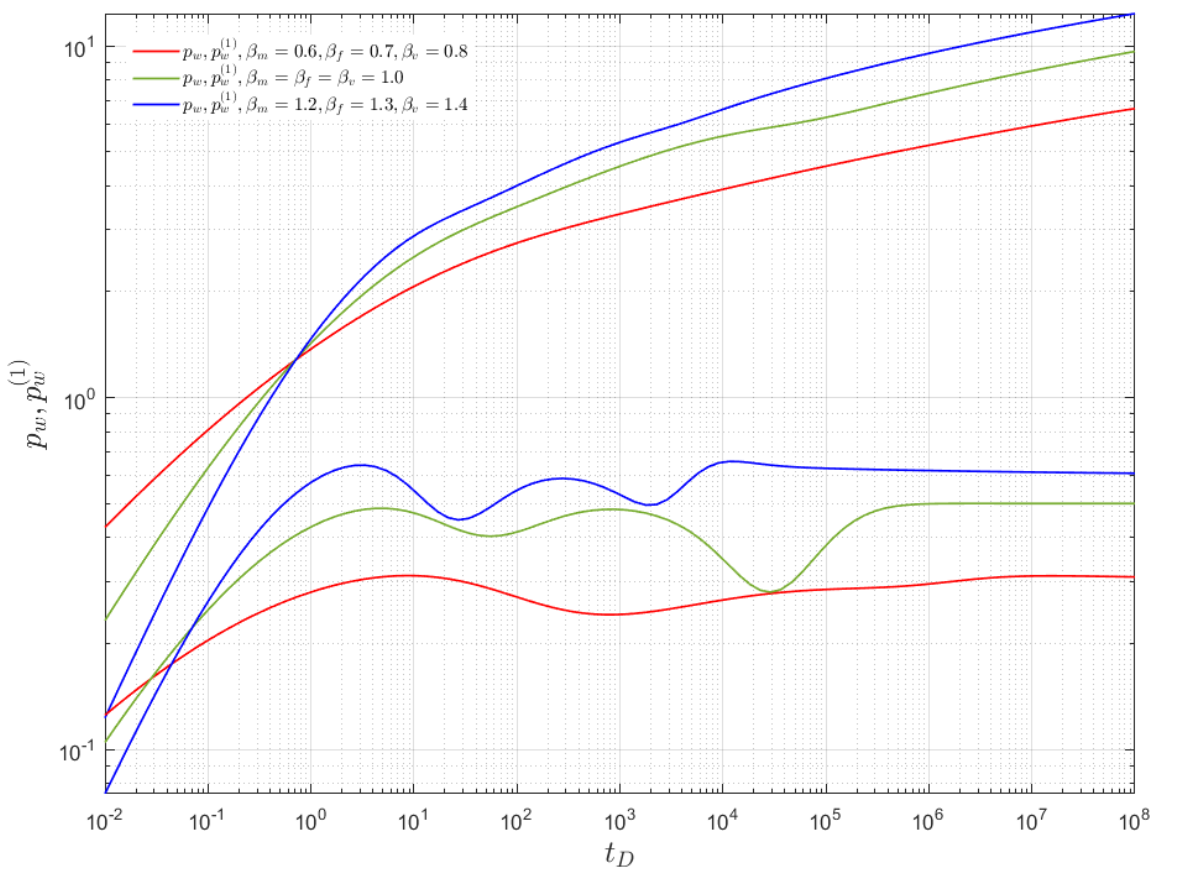}
     \end{subfigure}
        \caption{Behavior of pressure $p_w$ and pressure derivative $p^{(1)}_w$ for different values of $\beta_m, \ \beta_f$ and $\beta_v$.  All solutions were obtained with the particular values $\kappa_f=0.75,\ \kappa_v=\omega_f=0.02,\ \omega_v=0.8,\  \lambda_{mf}=10^{-3},\ \lambda_{vf}=10^{-8}$ and $ \lambda_{fv}=10^{-5}$. The solution in green in all the graphs corresponds to the classic case.}
        \label{fig:001}
\end{figure}

It should be mentioned that the numerical computation of an inverse Laplace transform is not easily done. In general, it is called a badly conditioned or ill-posed problem 
\cite{graf2012applied}. Some methods work quite well for certain functions, while for other functions they give poor results. As a general rule, numerical inversion methods work best for problems where the essential behavior of the original function is concentrated in a finite interval $[0, T]$, i.e. the original function should be decaying for increasing values of its independent variable \cite{graf2012applied}. Considering the aforementioned and the fact that the pressure $p_w$ is expected to satisfy the following condition:

\begin{eqnarray}
\lim_{r_d\to \infty}p_w \to 0,
\end{eqnarray}

then,  it is possible to implement the Stehfest algorithm    \cite{stehfest1970algorithm} in the equation \eqref{eq:001}. Finally, by obtaining the pressure $p_w$ in real space, it is possible to numerically approximate its derivative $p_w^{(1)}$ to obtain the graphs presented in Figure \ref{fig:001}.

\section{Conclusions}
%We have presented the methodology of a new model for pressure deficit in an oil will using derivatives Caputo  in the flow equation, through the appropriate mathematical machinery, to approach models that reflect the anomalous behavior of fluids in porous media, would allow to reflect the fractality of the medium or the media in the behavior of the fluid, a requirement would be to have the dimensions of the media. Some results about are from Casar-Gonzalez in \cite{Casar2000} and Hewett in \cite{Hewett1994}.
The memory formalism expressed through a generalized Darcy's law with Caputo derivative can facilitate reservoir modeling without the need to know their geometric structure to approach models that reflect anomalous behavior of fluids in porous media. The property of symmetry is important, because it allows reducing the problem to one dimension.The simplification of the equations to an algebraic problem facilitate the solution.

Partially funded by PAPPIT$\-$IT$\_$101421, UNAM.

\bibliography{biblio}

\begin{thebibliography}{10}

\bibitem{Warren1963}
J.E. Warren and P.J. Root.
\newblock The behavior of naturally fractured reservoirs.
\newblock {\em Society of Petroleum Engineers Journal}, 3(03):245--255, sep
  1963.

\bibitem{Chang1990}
J.~Chang and Y.~C. Yortsos.
\newblock Pressure transient analysis of fractal reservoirs.
\newblock {\em {SPE} Formation Evaluation}, 5(01):31--38, mar 1990.

\bibitem{CamachoVelazquez2008}
R.~Camacho~Velazquez, G.~Fuentes-Cruz, and M.~A. Vasquez-Cruz.
\newblock Decline-curve analysis of fractured reservoirs with fractal geometry.
\newblock {\em {SPE} Reservoir Evaluation {\&} Engineering}, 11(03):606--619,
  jun 2008.

\bibitem{Metzler1994}
R.~Metzler, W.~G. Gl\"{o}ckle, and T.~F. Nonnenmacher.
\newblock Fractional model equation for anomalous diffusion.
\newblock {\em Physica A: Statistical Mechanics and its Applications},
  211(1):13--24, oct 1994.

\bibitem{CamachoVelazquez2005}
R.~Camacho~Velazquez, M.~A. Vasquez-Cruz, R.~Castrejon-Aivar, and
  V.~Arana-Ortiz.
\newblock Pressure transient and decline curve behaviors in naturally fractured
  vuggy carbonate reservoirs.
\newblock {\em {SPE} Reservoir Evaluation {\&} Engineering}, 8(02):95--112, apr
  2005.

\bibitem{Ertekin1984}
T.~Ertekin.
\newblock Principles of numerical simulation of oil reservoirs {\textemdash} an
  overview.
\newblock In {\em Heavy Crude Oil Recovery}. Springer Netherlands, 1984.

\bibitem{martinez2017applications2}
B.~Mart{\'\i}nez-Salgado, R.~Rosas-Sampayo, A.~Torres-Hern{\'a}ndez, and
  C.~Fuentes.
\newblock Application of fractional calculus to oil industry.
\newblock {\em Fractal Analysis: Applications in Physics, Engineering and
  Technology}, 2017.
\newblock
  \url{https://www.intechopen.com/books/fractal-analysis-applications-in-physics-engineering-and-technology}.

\bibitem{fernando2017fractional}
F.~Brambila-Paz and A.~Torres-Hernandez.
\newblock Fractional newton-raphson method.
\newblock {\em arXiv preprint arXiv:1710.07634}, 2017.
\newblock \url{https://arxiv.org/pdf/1710.07634.pdf}.

\bibitem{brambila2018fractional}
F.~Brambila-Paz, A.~Torres-Hernandez, U.~Iturrar{\'a}n-Viveros, and
  R.~Caballero-Cruz.
\newblock Fractional newton-raphson method accelerated with aitken's method.
\newblock {\em arXiv preprint arXiv:1804.08445}, 2018.
\newblock \url{https://arxiv.org/pdf/1804.08445.pdf}.

\bibitem{torreshern2020}
A.~Torres-Hernandez, F.~Brambila-Paz, and E.~{De-la}-Vega.
\newblock Fractional newton-raphson method and some variants for the solution
  of nonlinear systems.
\newblock {\em Applied Mathematics and Sciences: An International Journal
  (MathSJ)}, 2020.
\newblock \url{https://airccse.com/mathsj/papers/7120mathsj02.pdf}.

\bibitem{torres2020reduction}
A.~Torres-Hernandez, F.~Brambila-Paz, and P.M. Rodrigo.
\newblock Reduction of a nonlinear system and its numerical solution using a
  fractional iterative method.
\newblock {\em arXiv preprint arXiv:2007.02776}, 2020.
\newblock \url{https://arxiv.org/pdf/2007.02776.pdf}.

\bibitem{torres2020fractional}
A.~Torres-Hernandez, F.~Brambila-Paz, P.M. Rodrigo, and E.~De-la Vega.
\newblock Fractional pseudo-newton method and its use in the solution of a
  nonlinear system that allows the construction of a hybrid solar receiver.
\newblock {\em Applied Mathematics and Sciences: An International Journal
  (MathSJ)}, 2020.
\newblock \url{https://airccse.com/mathsj/papers/7220mathsj01.pdf}.

\bibitem{martineznumerical}
C.~A. Mart{\i}nez and F.~Brambila-Paz.
\newblock Numerical comparison between rbf schemes with respect to other
  approaches to solve fractional partial differential equations and their
  advantages when choosing non-uniform nodes.
\newblock {\em Journal of Mathematics and Statistical Science}, 5:85--105, May
  2019.

\bibitem{martinez2017applications1}
C.~A. Mart{\'\i}nez and C.~Fuentes.
\newblock Applications of radial basis function schemes to fractional partial
  differential equations.
\newblock {\em Fractal Analysis: Applications in Physics, Engineering and
  Technology}, 2017.
\newblock
  \url{https://www.intechopen.com/books/fractal-analysis-applications-in-physics-engineering-and-technology}.

\bibitem{LeMehaute1984}
A.~Le~Mehaute.
\newblock Transfer processes in fractal media.
\newblock {\em Journal of Statistical Physics}, 36(5-6):665--676, sep 1984.

\bibitem{Raghavan2013}
R.~Raghavan and C.~Chen.
\newblock Fractional diffusion in rocks produced by horizontal wells with
  multiple, transverse hydraulic fractures of finite conductivity.
\newblock {\em Journal of Petroleum Science and Engineering}, 109:133--143,
  2013.

\bibitem{Holy2016}
R.~W: Holy.
\newblock {\em Numerical investigation of 1D anomalous diffusion in fractured
  nanoporous reservoirs}.
\newblock {PhD} {D}issertation, Colorado School of Mines. Arthur Lakes Library,
  2016.

\bibitem{hatano2013determination}
Y.~Hatano, J.~Nakagawa, S.~Wang, and M.~Yamamoto.
\newblock Determination of order in fractional diffusion equation.
\newblock {\em J. Math-for-Ind. A}, 5(51):118, 2013.

\bibitem{Bear1988}
J.~Bear.
\newblock {\em Dynamics of fluids in porous media}.
\newblock Dover, 1988.

\bibitem{Peaceman1977}
D.~W. Peaceman.
\newblock {\em Fundamentals of Numerical Reservoir Simulation}, volume~6 of
  {\em Developments in Petroleum Science}.
\newblock Elsevier Science, 1st edition, 1977.

\bibitem{Fuentes2014}
C.~Fuentes.
\newblock {Informe VI. Fondo Sectorial CONACYT-SENER-Hidrocarburos}.
\newblock Technical Report S0018-2011-11, jun 2014.

\bibitem{raghavan2018conceptual}
Rajagopal Raghavan, Chih-Cheng Chen, et~al.
\newblock A conceptual structure to evaluate wells producing fractured rocks of
  the permian basin.
\newblock In {\em SPE Annual Technical Conference and Exhibition}. Society of
  Petroleum Engineers, 2018.

\bibitem{Baleanu2011}
D.~Baleanu, E.~Diethelm, K.~Scalas, and J.J. Trujillo.
\newblock {\em Fractional Calculus}.
\newblock WORLD SCIENTIFIC, 2012.

\bibitem{graf2012applied}
Urs Graf.
\newblock {\em Applied Laplace transforms and z-transforms for scientists and
  engineers: a computational approach using a Mathematica package}.
\newblock Birkh{\"a}user, 2012.

\bibitem{stehfest1970algorithm}
Harald Stehfest.
\newblock Algorithm 368: Numerical inversion of laplace transforms [d5].
\newblock {\em Communications of the ACM}, 13(1):47--49, 1970.

\end{thebibliography}
\bibliographystyle{unsrt}

\end{document}